\documentclass[11pt]{article}
\usepackage[utf8]{inputenc}
\usepackage[a4paper, margin = 0.5in, bottom = 1in]{geometry}
\usepackage{soul}
\usepackage{authblk}

\date{}
\usepackage{amsmath}
\usepackage{amsfonts}
\usepackage{amssymb}
\usepackage{mathrsfs}
\usepackage{bm}
\usepackage{multicol}
\usepackage{wasysym}
\usepackage{graphicx}
\usepackage{caption}
\usepackage{cite}
\usepackage{subcaption}
\usepackage{skak}
\usepackage[table]{xcolor}
\usepackage{bbm}

\usepackage{amsfonts}
\usepackage{subcaption}
\usepackage{lmodern}
\usepackage{xcolor}
\usepackage{csquotes}
\usepackage{mathtools}
\usepackage{float}
\usepackage{amsmath,amssymb,amsthm,comment}
\usepackage{amsfonts}
\usepackage{subcaption}
\usepackage{xcolor}
\usepackage{csquotes}
\usepackage{mathtools}
\usepackage{float}
\usepackage{amsmath,amssymb,amsthm}
\usepackage[colorlinks,citecolor=red,urlcolor=blue,linkcolor=blue]{hyperref}
\usepackage{mathrsfs,dsfont}
\usepackage{graphicx}
\usepackage{enumerate}
\usepackage{nicefrac}
\usepackage[font={normalsize},width=1.12\linewidth]{caption}
\usepackage{comment}
\usepackage{hyperref}
\usepackage{pgfplots}
\pgfplotsset{compat=1.15}
\usepackage{mathrsfs}
\usetikzlibrary{arrows}

\newtheorem{nummer}{ }
\newtheorem{thm}[nummer]{Theorem}
\newtheorem{prp}[nummer]{Proposition}
\newtheorem{lem}[nummer]{Lemma}

\newtheorem{fct}[nummer]{Fact}
\newtheorem{q}[nummer]{Question}
\newtheorem{defi}[nummer]{Definition}

\DeclareMathOperator{\vvec}{vec}

\makeatletter
\def\opargproof[#1]{\par\noindent {\bf #1 }}
 \makeatother

\begin{document}
\medskip\medskip
\begin{center}
\vspace*{50pt}
{\LARGE\bf Neighbour Sum Patterns : Chessboards to Toroidal Worlds}

\bigskip
{\small Sayan Dutta}\\[1.2ex] 
{\scriptsize Department of Mathematics and Statistics, Institute of Science Education and Research (IISER) Kolkata\\
Département de mathématiques et de statistique, Université de Montréal\\
sayan.dutta@umontreal.ca\\
\href{https://sites.google.com/view/sayan-dutta-homepage}{https://sites.google.com/view/sayan-dutta-homepage}}
\\[1.8ex]

{\small Ayanava Mandal}\\[1.2ex] 
{\scriptsize School of Mathematical Sciences, NISER, Jatni, Odisha 752050, India\\ ayanava.mandal@niser.ac.in}\\[1.8ex]

{\small Sohom Gupta}\\[1.2ex] 
{\scriptsize Department of Atmospheric and Oceanic Sciences, McGill University.\\
Department of Physical Sciences, Institute of Science Education and Research (IISER) Kolkata\\
sohom.gupta@mail.mcgill.ca}\\[1.8ex]

{\small Sourin Chatterjee}\\[1.2ex] 
{\scriptsize Department of Mathematics and Statistics, Institute of Science Education and Research (IISER) Kolkata\\
Institut de Neurosciences des Systèmes (INS), UMR1106, Aix-Marseille Université, Marseilles, France\\
sourin.chatterjee@univ-amu.fr}
\\[1.8ex]

\end{center}

\hspace{5ex}{\small{\it key-words\/}: p-adic valuation, Kronecker product, spatial lattice, discrete harmonic function, cyclotomic polynomials}

\begin{abstract}
We say that a chessboard filled with integer entries is a \textit{solution} to the \textit{neighbor-sum problem} if the number appearing on each cell is the sum of entries in its neighboring cells, where neighbors are cells sharing a common edge or vertex. We show that an $n\times n$ chessboard satisfies this property if and only if $n\equiv 5\pmod 6$. We investigate the existence of solutions of rectangular, toroidal boards, as well as on Neumann neighborhoods, including a nice connection to discrete harmonic functions. We construct solutions on infinite boards are also presented. Finally, awe explore three-dimensional analogs of these boards using properties of cyclotomic polynomials and relevant ideas conjectured.  
\end{abstract}

\section{Introduction.}
The Regional Mathematical Olympiad (RMO) is the second of a series of math tests held in India, which all leads up to participation in the International Mathematical Olympiad (IMO). 

Our inspiration stems from \href{https://www.isical.ac.in/~rmo/papers/rmo/rmo-1991.pdf}{\textcolor{blue}{RMO 1991}} Problem 8 -

\textit{The 64 squares of an 8 × 8 chessboard are filled with positive integers in such a way that each
integer is the average of the integers on the neighboring squares. (Two squares are neighbors
if they share a common edge or a common vertex. Thus a square can have 8, 5 or 3 neighbors
depending on its position). Show that all the 64 integer entries are in fact equal.}

\textbf{Brief Solution}: \textit{Any given entry must lie in between the smallest and largest entries of its neighbors.
Thus, the largest entry on the board must be surrounded by identical entries.
This forces all entries to be equal.}

While this has a surprisingly mundane answer, a small modification to the criterion might not be so! This is the analog we explore in this paper:  

\textit{An $n \times n$ chessboard for $n \geq 3$, with each square bearing an integer, is said to be a \textbf{solution} to the neighbor sum problem if each number is the sum of the numbers on the neighboring squares. Two squares are neighbors if they share a common edge or a common vertex.\footnote{This definition of a \textit{neighborhood} is referred to as the Moore neighborhood. Often in network analysis and cellular automata problems, people consider neighborhoods to be restricted to squares that share an edge\textemdash which is called a Neumann neighborhood and seen in the context of the problem later.} A chessboard with all squares bearing the number zero is said to be a \textbf{trivial} solution. How many such non-trivial solutions exist for a chessboard of given size $n\times n$?}

At first glance, the problem might seem to be rooted in combinatorics, but certain observations favor a different angle. It is clear that given a non-trivial solution $((x_{ij}))$ of the matrix representation $X$ of a chessboard (of dimension $n\times n$ for some $n$), all matrix representations of the form $\alpha X=((\alpha x_{ij})),\;\alpha \in\mathbb{Z}\backslash\{0\}$ are also valid non-trivial solutions. Furthermore, the sum of two solutions is also another solution. This motivates the idea of a transformation which contains the vectorization \cite{Dhrymes2000} of $X$ in its kernel.\footnote{The vectorization of a matrix $A$, denoted by $\bm{\vvec (A)}$, is a vector formed by stacking the columns of A in a top-down format.}

Section 2 documents an attempt to figure out all solutions for square chessboards of size $n\times n$. This combines the classical way to study eigenvalues of adjacency matrices of graphs with some standard results on the Kronecker product of matrices (Facts \ref{fct:4},\ref{fct:5}) and some $p$-adic analysis (Lemma \ref{speyer}). These ideas are then reformulated appropriately for certain generalizations, such as rectangular, toroidal, infinite 2D, and higher dimensional chessboards which are solutions to the neighbor sum problem. Those span Sections 3 to 9.

\begin{figure}[htp!]
    \centering
    \includegraphics[width=0.6\textwidth]{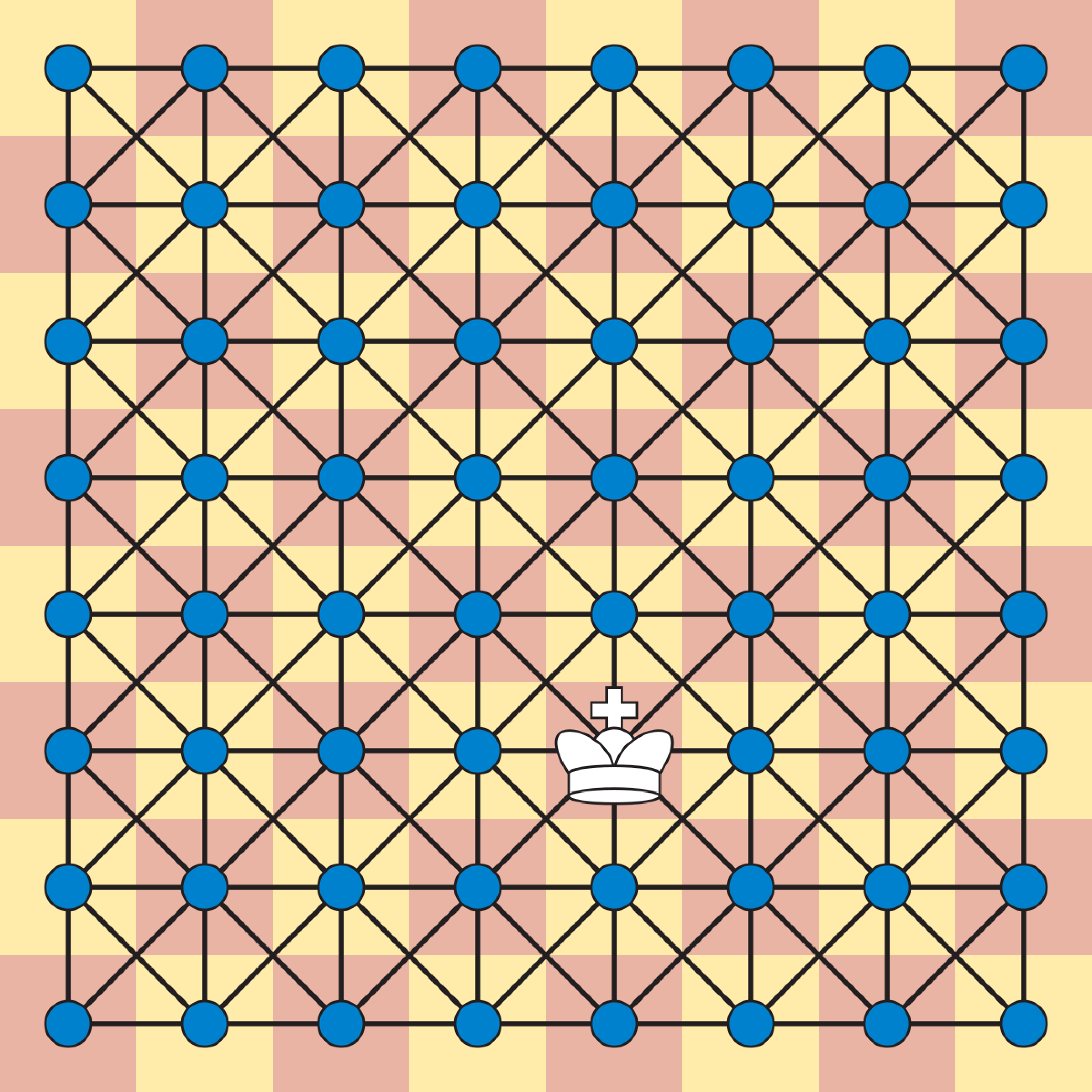}
    \caption{King's Graph on a standard $8\times 8$ chessboard. Image Courtesy: \href{https://commons.wikimedia.org/wiki/User:David_Eppstein}{David Epstein.}}
\end{figure}

\section{Finding square boards with such solutions}
Trying to resolve our problem has led us to a critical observation - the sum of two solutions is a solution itself. It remains to produce an appropriate transformation $T$ such that for a chessboard $\symking\in M_n(\mathbb{Z})$ (the use of this notation being a subtle foreshadowing) with its vectorization $\bm{\vvec(\symking)}\in \mathbb{Z}^{n^2}$, one has
\[T\bm{\vvec(\symking)} = \bm{0}\]

An immediate transformation is one that replaces every element with the sum of its neighbors minus the element itself. Clearly, all solutions would be in the kernel of this transformation. For a $2\times 2$ chessboard,
\[\def\arraystretch{1.5}
        \symking_{2\times 2} = \begin{array}{|c|c|}\hline
            x_1 & x_3 \\\hline
            x_2 & x_4 \\\hline
        \end{array}
\]
and the corresponding $4\times 4$ transformation would be given by
\[T_{4\times 4} = \begin{pmatrix}
    -1&1&1&1\\
    1&-1&1&1\\
    1&1&-1&1\\
    1&1&1&-1
\end{pmatrix} .\]

This is non-singular, so this does not have a non-trivial kernel. Note that we are only interested in $n\geq 3$, and the same can be shown for $T_{3\times 3}$ for $n=3$. But are there some values of $n$ for which non-trivial solutions can be present? This takes us to our main result.

\subsection{Existence of solutions.}

\begin{thm}[\textbf{Existence}]\label{thm:square}
    An $n\times n$ chessboard has a non-trivial solution of the neighbor-sum property if and only if $6 \mid (n + 1)$.
\end{thm}

This is the key result for this section, and its proof requires certain machinery that is introduced and discussed as follows.

Let $\symking$ be an $n \times n$ chessboard, and let $T_n$ be the transformation mentioned above. Then, $A_n = T_n + \mathbb{I}_{n^2}$ is precisely the adjacency matrix of the square chessboard where adjacency is only amongst neighbors. Denote the $i$-th element of $\bm{\vvec(\symking)}$ as $\symking_i$. Then the $i$-th entry of the vector $A_n\bm{\vvec{\symking}}$ gives the sum of the neighbors of the square $\symking_i$. This means that the neighbor-sum property can be expressed as $A_n\bm{\vvec\symking} = \bm{\vvec{\symking}}$.

\begin{prp}\label{prp:solutions_kernel}
    {The set of all vectorised solutions is precisely $\ker(T_n)$, which always contains the trivial solution.}
\end{prp}

\begin{prp}\label{prp:kings_graph}
    Define $B_n \in M_n(\mathbb{Z})$ with $B_{ij} = 1$ when $|i - j| \leq 1$ and $B_{ij} = 0$ otherwise. Then, \[
        A_n = B_n\otimes B_n -  \mathbb{I}_{n^2}, \qquad
        T_n = B_n\otimes B_n - 2\mathbb{I}_{n^2}.
    \]
\end{prp}
\begin{proof}
    Note that $B_n$ can be interpreted as the adjacency matrix (\cite{biggs1993algebraic}, pg. 7) of a graph $G_n$ on vertices $\{1, \dots, n\}$, where $i, i'$ are neighbors when $|i - i'| \leq 1$. In the Cartesian product (\cite{hahn2013graph}, pg. 115 - 116) $G_n\times G_n$, whose adjacency matrix is $B_n\otimes B_n$, we have an edge between $(i, j)$ and $(i', j')$ precisely when $|i - i'| \leq 1$ and $|j - j'| \leq 1$. Removing self-loops by subtracting $\mathbb{I}_{n^2}$ from the adjacency matrix yields the transformation $A_n$ as desired.
\end{proof}
\textit{Remark.} The graph $G_{n,n}$ is called the King's Graph \cite{chang1998algorithmic,berend2005two} because it shows the movement of a King on a chessboard. At any square (equivalent to a node in the graph), the King has 3, 5, or 8 adjacent squares to move to depending on its position on the board. The adjacency matrix of this graph is $A_n=B_n\otimes B_n - \mathbb{I}_{n^2}$. This implicit relation motivates the $\symking$ notation.

\[
        B_5 = \begin{bmatrix}
            1 & 1 & 0 & 0 & 0 \\
            1 & 1 & 1 & 0 & 0 \\
            0 & 1 & 1 & 1 & 0 \\
            0 & 0 & 1 & 1 & 1 \\
            0 & 0 & 0 & 1 & 1
        \end{bmatrix},
        \qquad
        B_5\otimes B_5 = \begin{bmatrix}
            B_5 & B_5 & 0 & 0 & 0 \\
            B_5 & B_5 & B_5 & 0 & 0 \\
            0 & B_5 & B_5 & B_5 & 0 \\
            0 & 0 & B_5 & B_5 & B_5 \\
            0 & 0 & 0 & B_5 & B_5
        \end{bmatrix}.
\]
    
With this, our search for non-trivial chessboards with the neighbor-sum property reduces to finding eigenvectors of $B_n\otimes B_n$ corresponding to the eigenvalue $2$.

\begin{fct}[Henderson, Pukelsheim and Searle \cite{henderson1983history}]\label{fct:4}
    The eigenvalues of $A\otimes B$ are $\{\lambda_i\mu_j\}$, where $\{\lambda_i\}$ are the eigenvalues of $A$, and $\{\mu_j\}$ are the eigenvalues of $B$.
\end{fct}

\begin{fct}[Kulkarni, Schmidt and Tsui \cite{kulkarni1999eigenvalues}] \label{fct:5}
    The eigenvalues of $B_n$ are $\lambda_k = 1 + 2\cos(k\pi / (n + 1))$ for $k = 1, \dots, n$. This is due to the tridiagonal Toeplitz form of $B_n$ with all non-zero elements being unity.
\end{fct}

Using these two facts, we formulate the following proposition.

\begin{prp}\label{prp:eigenvalue_eq}
    The space $\ker(T_n)$ is non-trivial if and only if there exist $p, q \in \mathbb{N}$ such that \[
        \left(1 + 2\cos\left(\frac{p\pi}{n + 1}\right)\right) \left(1 + 2\cos\left(\frac{q\pi}{n + 1}\right)\right) = 2
    \]
    and $1 \leq p, q \leq n$.
\end{prp}

In order to deal with the equation in Proposition \ref{prp:eigenvalue_eq} and others similar to it, we require the following result.

\begin{thm}\label{thm:rational_solutions}
    The only solutions of \[
        \left(1 + 2\cos(u\pi)\right) \left(1 + 2\cos(v\pi)\right) = 2.
    \] where $u, v \in \mathbb{Q} \cap (0, 1)$ are $u = 1/3, v = 1/2$ and $u = 1/2, v = 1/3$.
\end{thm}

\begin{proof}
    The given equation can be rewritten as
    \begin{equation}\label{product}
        (\alpha+1+\alpha^{-1}) (\beta + 1 + \beta^{-1}) = 2,
    \end{equation} where $\alpha = e^{iu\pi}, \beta = e^{iv\pi}$ are roots of unity with positive imaginary parts.

    Let $u = p / N$, $v = q / N$ be a solution of this equation, where $p, q, N \in \mathbb{N}$ and $1 \leq p, q < N$. Set $R = \mathbb{Z}[e^{\pi i/N}]$. Let $\mathfrak{p}$ be a prime of $R$ lying over the prime $2$ in $\mathbb{Z}$, and let $v : R \to \mathbb{Q}$ be the $\mathfrak{p}$-adic valuation (\cite{dummit2004abstract} pg. 755), normalized so that $v(2)=1$. So, Equation~\ref{product} gives
    \begin{equation}\label{sum}
        v(\alpha+1+\alpha^{-1}) + v(\beta+1+\beta^{-1}) = 1.
    \end{equation}

    To proceed any further, we state the following lemma.
    
    \begin{lem}\label{speyer}
        Let $\eta$ be a primitive $m$-th root of unity. Then 
    $$v(\eta+1+\eta^{-1}) = \begin{cases} \infty & m = 3 \\ 1/2^k & m = 3 \cdot 2^{k+1} \ \text{for} \ k \geq 0 \\ 0 & \text{otherwise} \end{cases}.$$
    \end{lem}

    Now, it is clear that the only ways to decompose $1$ as a sum of two numbers in $\{ \infty, 1, 1/2, 1/4, \cdots, 0 \}$ are $1+0$, $0+1$,  and $1/2 + 1/2$. We use this in Equation~\ref{sum}.

    \textbf{Case I}: If $v(\alpha + 1 + \alpha^{-1}) = 1$, $v(\beta + 1 + \beta^{-1}) = 0$, then $\alpha$ must be a primitive $6$-th root of unity, forcing $\alpha = e^{i\pi / 3}$, $u = 1 / 3$. This in turn forces $v = 1 / 2$. Interchanging the roles of $\alpha, \beta$ yields the solution $u = 1 / 2$, $v = 1 / 3$.

    \textbf{Case II}: If $v(\alpha + 1 + \alpha^{-1}) = v(\beta + 1 + \beta^{-1}) = 1 / 2$, then $\alpha, \beta$ must be primitive $12$-th roots of unity, forcing $\alpha, \beta \in \{e^{i\pi / 6}, e^{5i\pi / 6}\}$. This gives, $\alpha + 1 + \alpha^{-1}, \beta + 1 + \beta^{-1} \in \{1 + \sqrt{3}, 1 - \sqrt{3}\}$; but these do not satisfy Equation~\ref{product}.
\end{proof}

Now, we return to the proof of the Lemma we just used.\footnote{Note that this proof uses the language of Algebraic Number Theory. We sincerely believe that it is possible to formulate a different proof using only elementary tools, but we will leave that to the enthusiastic readers.}
\begin{proof}[Proof of Lemma~\ref{speyer}]
    The case $\eta=1$ (so $m=1$) is easy to check by hand, so we assume that $\eta \neq 1$ from now on.

    We have
    $$\eta+1+\eta^{-1} = \eta^{-1}\cdot \frac{\eta^3-1}{\eta-1}$$
    so
    $$v(\eta+1+\eta^{-1}) = v(\eta^3-1) - v(\eta-1).$$
    Now, if $\omega$ is a primitive $\ell$-th root of unity, then
    $$v(\omega-1) = \begin{cases} \infty & \ell=1 \\ 1/2^k & \ell = 2^{k+1} \\ 0 & \text{otherwise} \end{cases}$$
    as, if $ \ell=1 $, then $ \zeta_\ell - 1 =0$; if two distinct primes divide $ \ell $, then $ 1-\zeta_\ell $ is a unit; and if $\ell=p^k $, then we have $ (\zeta_\ell - 1)^{\phi(\ell)} = (p) $ as ideal and hence, $ v(\zeta_\ell -1) $ is nonzero only when $ p=2 $ and the valuation is $ \frac{1}{\phi(\ell)} $.

    Since $\eta$ is a primitive $m$-th root of unity, $\eta^3$ is a primitive $\frac{m}{\gcd (m,3)}$-th root of unity, so combining the above two equations gives the claim.
\end{proof}

{Finally, we are equipped with enough tools to provide a proof for Theorem 1.}

\begin{proof}[Proof of Theorem~\ref{thm:square}]
    Propositions \ref{prp:solutions_kernel} and \ref{prp:eigenvalue_eq} give a criterion for the existence of non-trivial $n\times n$ chessboards with the neighbor-sum property.
    Theorem \ref{thm:rational_solutions} forces both $2\mid (n + 1)$ and $3\mid(n + 1)$, from which the claim follows.
\end{proof}

\subsection{Looking at the solutions.}

If the dimension $n$ of a square board is one less than a multiple of 6, the transformation $B_n\otimes B_n$ has eigenvalue $2$ with multiplicity $2$, which implies that $\ker(T_n)$ is two-dimensional. It is a matter of computation to yield the exact solutions, which correspond to the associated eigenvectors. For the smallest board satisfying the neighbor-sum property ($n=5$),

\[\def\arraystretch{1.7}\arraycolsep=1pt
        \symking^{(1)}_{5\times 5} =
        \begin{array}{|c|c|c|c|c|}\hline
            1 & 0 & -1 & 0 & 1 \\\hline
            1 & 0 & -1 & 0 & 1 \\\hline
            \phantom{-}0\phantom{-} & \phantom{-}0\phantom{-} & \phantom{-}0\phantom{-} & \phantom{-}0\phantom{-} & \phantom{-}0\phantom{-} \\\hline
            -1 & 0 & 1 & 0 & -1 \\\hline
            -1 & 0 & 1 & 0 & -1 \\\hline
        \end{array}
        \text{ and }
        \def\arraystretch{1.7}\arraycolsep=1pt
        \symking^{(2)}_{5\times 5} =
        \begin{array}{|c|c|c|c|c|}\hline
            1 & 1 & 0 & -1 & -1 \\\hline
            \phantom{-}0\phantom{-} & \phantom{-}0\phantom{-} & \phantom{-}0\phantom{-} & \phantom{-}0\phantom{-} & \phantom{-}0\phantom{-} \\\hline
            -1 & -1 & 0 & 1 & 1 \\\hline
            0 & 0 & 0 & 0 & 0 \\\hline
            1 & 1 & 0 & -1 & -1 \\\hline
        \end{array}
    \]

are the only solutions.

Clearly, these solutions are the transposes of each other, and have distinct vectorizations which makes them distinct solutions. These are surprisingly simple solutions, with every element being in the set $\{-1,0,1\}$. Custom solutions can be produced with linear combinations $\lambda \symking^{(1)}_{5\times 5} + \mu\symking^{(2)}_{5\times 5}$ for $\lambda,\mu \in \mathbb{Z}$, both not zero.

One observation is critical in understanding the kind of solutions one should expect to see for larger square boards admitting solutions. Define a \emph{phantom boundary} to be a boundary of cells (all containing zero) of a board such that the elements of the boundary contribute only to the neighbor-sum property of the board and not themselves.\footnote{For a $p\times p$ board, the boundary is the collection of rows $1,p$ and columns $1,p$.}

Consider the following $2\times 2$ board with a phantom boundary represented in white around the yellow (gray in print) board, yielding a $4\times 4$ board.

\[\def\arraystretch{1.5}
        \begin{array}{|c|c|c|c|}\hline
            0&0&0&0\\\hline
            0&\cellcolor[HTML]{F7F978}x_1 & \cellcolor[HTML]{F7F978}x_3&0 \\\hline
            0&\cellcolor[HTML]{F7F978}x_2 & \cellcolor[HTML]{F7F978}x_4&0 \\\hline
            0&0&0&0\\\hline
        \end{array}
\]

The presence of the phantom boundary does not alter the conditions necessary for the $2\times 2$ to satisfy the neighbor-sum property as the zeroes don't contribute to the sum. This idea plays a critical role in identifying disjoint solutions in large boards.

We can already see that the $5\times 5$ solutions can be formed of small $2\times 1$ units emphasized in the following figure.

\[\def\arraystretch{1.7}\arraycolsep=1pt
        \symking^{(1)}_{5\times 5} =
        \begin{array}{|c|c|c|c|c|}\hline
            \cellcolor[HTML]{F7F978}1 & 0 & \cellcolor[HTML]{F7F978}-1 & 0 & \cellcolor[HTML]{F7F978}1 \\\hline
            \cellcolor[HTML]{F7F978}1 & 0 & \cellcolor[HTML]{F7F978}-1 & 0 & \cellcolor[HTML]{F7F978}1 \\\hline
            \phantom{-}0\phantom{-} & \phantom{-}0\phantom{-} & \phantom{-}0\phantom{-} & \phantom{-}0\phantom{-} & \phantom{-}0\phantom{-} \\\hline
            \cellcolor[HTML]{F7F978}-1 & 0 & \cellcolor[HTML]{F7F978}1 & 0 & \cellcolor[HTML]{F7F978}-1 \\\hline
            \cellcolor[HTML]{F7F978}-1 & 0 & \cellcolor[HTML]{F7F978}1 & 0 & \cellcolor[HTML]{F7F978}-1 \\\hline
        \end{array}
    \]
Adding a phantom boundary to this solution clearly shows that the solution can be divided into 6 disjoint regions (three $(1,1)$ and three $(-1,-1)$, alternating) separated by zeroes. This pattern can easily be repeated to get the two solutions for $n=11,17,\ldots$. Since the kernel is always 2 dimensional, the solutions on $n\times n$ formed from extensions of those on $5\times 5$ form a basis for the eigenspace. This gives us a complete characterization of solutions of the neighbor-sum property on square boards.

It is easy to see and prove that in the standard solutions, every second column and third row in $\symking_{n\times n}^{(1)}$ (respestively second row and third column in $\symking^{(2)}_{n\times n}$) contains only zero elements. The only zeroes that are common are at positions $(i,j)$ where either both $i$ and $j$ are even or they are both multiples of 3. So, any non-trivial linear combination $\lambda \symking^{(1)}_{5\times 5} + \mu\symking^{(2)}_{5\times 5}$ for $\lambda,\mu \in \mathbb{Z}$ would preserve the zeroes in those positions.

\textbf{\textit{Remark.}} As the standard solutions form a basis for the kernel, there cannot be any square board with a non-trivial solution without any zero elements.

\textbf{Overview of the Paper.} This previous remark answers our original question stated in \textbf{Section 1}. The techniques and results obtained in \textbf{Section 2}, especially Theorem \ref{thm:rational_solutions} and Lemma \ref{speyer} can be used to tackle the neighbor sum problems on non-square chessboards, even under different boundary conditions. \textbf{Section 3} deals with $m\times n$ rectangular chessboards and how their solutions can be constructed of the same $2\times 1$ solutions we noticed in the square board case. \textbf{Section 4} is a similar attempt at finding the existence criteria for solutions to toroidal boards.\\
\textbf{Section 5} looks at boards under the Neumann neighbor sum (Nns) property, where only edge adjacent squares are considered as neighbors. We state certain results for square boards and briefly study an interesting variation concerning harmonic functions on toroidal boards.\\
In \textbf{Sections 6} through \textbf{9}, we look at infinite two dimensional and finite higher dimensional analogs for the neighbor sum problem. While for an infinite chessboard, it is quite easy to generate solutions, higher dimensional analogs are far more complicated and require the use of advanced Algebraic results and ideas which are currently beyond our scope. We explore pathways through which this can be achieved, and remark on how the neighbor sum problem on any discrete graph can be studied.

\section{Reducing symmetry\textemdash from squares to rectangles}
A simple generalization of the neighbor-sum property can be made to rectangular boards of size $m\times n$ with $m,n\geq 2$.\footnote{For $m\neq1$, $n=1$, we get a one-dimensional strip, which has solutions when $m\equiv 2\;(\operatorname{mod} 3)$, with solutions easily constructible from the $2\times 1$ units on the square board.} Theorem~\ref{thm:rational_solutions} will play a role in finding solutions here.

Following similar arguments as in the case of $n\times n$ chessboards, it is not difficult to arrive at analogs of Propositions~\ref{prp:kings_graph} and \ref{prp:eigenvalue_eq} for $m \times n$ chessboards.

\begin{prp}
    The set of all $m\times n$ solutions i.e,, $m\times n$ chessboards with the neighbor-sum property, is precisely $\ker(T_{m,n})$, where \[
        T_{m, n} = B_m\otimes B_n - 2\mathbb{I}_{mn}.
    \]
\end{prp}

\begin{prp}\label{prp:ker_Tmn}
    The space $\ker(T_{m, n})$ is non-trivial if and only if there exist $p, q \in \mathbb{N}$ such that \[
        \left(1 + 2\cos\left(\frac{p\pi}{m + 1}\right)\right) \left(1 + 2\cos\left(\frac{q\pi}{n + 1}\right)\right) = 2.
    \]
    and $1 \leq p \leq m$, $1 \leq q \leq n$.
\end{prp}

This yields the following characterization.

\begin{thm}\label{thm:rectangle}
    Non-trivial $m\times n$ chessboards satisfying the neighbor-sum property exist if and only if $2 \mid (m+1)$ and $3 \mid (n+1)$.
\end{thm}
\begin{proof}
    Follows immediately from Proposition~\ref{prp:ker_Tmn} and Theorem~\ref{thm:rational_solutions}.
\end{proof}
\textit{\textbf{Remark.} The dimension of $\ker(T_{m, n})$ is at most $2$. It is equal to $1$ only if $m\neq n$.}

This remark follows directly from the solution space explored before. There are two fundamental solutions for a square board, and they can only \textit{fit} in a rectangular board if both dimensions are large enough. In all other cases, we have the nullity to be at most $1$. 

Some simple consequences of these are:
\begin{itemize}
    \item[1.] If $\symking_{m\times n}$ has a non-trivial kernel of dimension $d\leq 2$, then $\symking_{n\times m}$ also has a non-trivial kernel of dimension $d$. Further, the solutions are transposes of one-another.
    \item[2.] The standard solutions of a square board can be partitioned into disjoint non-trivial rectangular solutions.
    \item[3.] A chessboard of dimensions $m\times n$, where $m+1$ is even and $n+1$ is an odd multiple of 3 (or vice versa), has solution(s) by Theorem~\ref{thm:rectangle}. Then, a board of dimensions $(m+1)\times (n+1)$ also has solution(s). Furthermore, if a standard solution of the $m\times n$ board is made of $2\times 1$ units, then a \emph{corresponding} standard solution on the $(m+1)\times (n+1)$ board is made of similar $1\times 2$ units. This correspondence is clear when $\ker(T_{m,n})$ has dimension $1$.
\end{itemize}

\section{No boundaries now\textemdash Looking at Torii.}
For any finite board, boundary conditions are important to construct solutions, where the \textit{phantom-boundary} visualisation has come handy. This prompts a natural curiosity - what if there were no boundaries? And where else to look for but a good old torus!

\begin{figure}[htp!]
    \centering
    \includegraphics[width=0.6\textwidth]{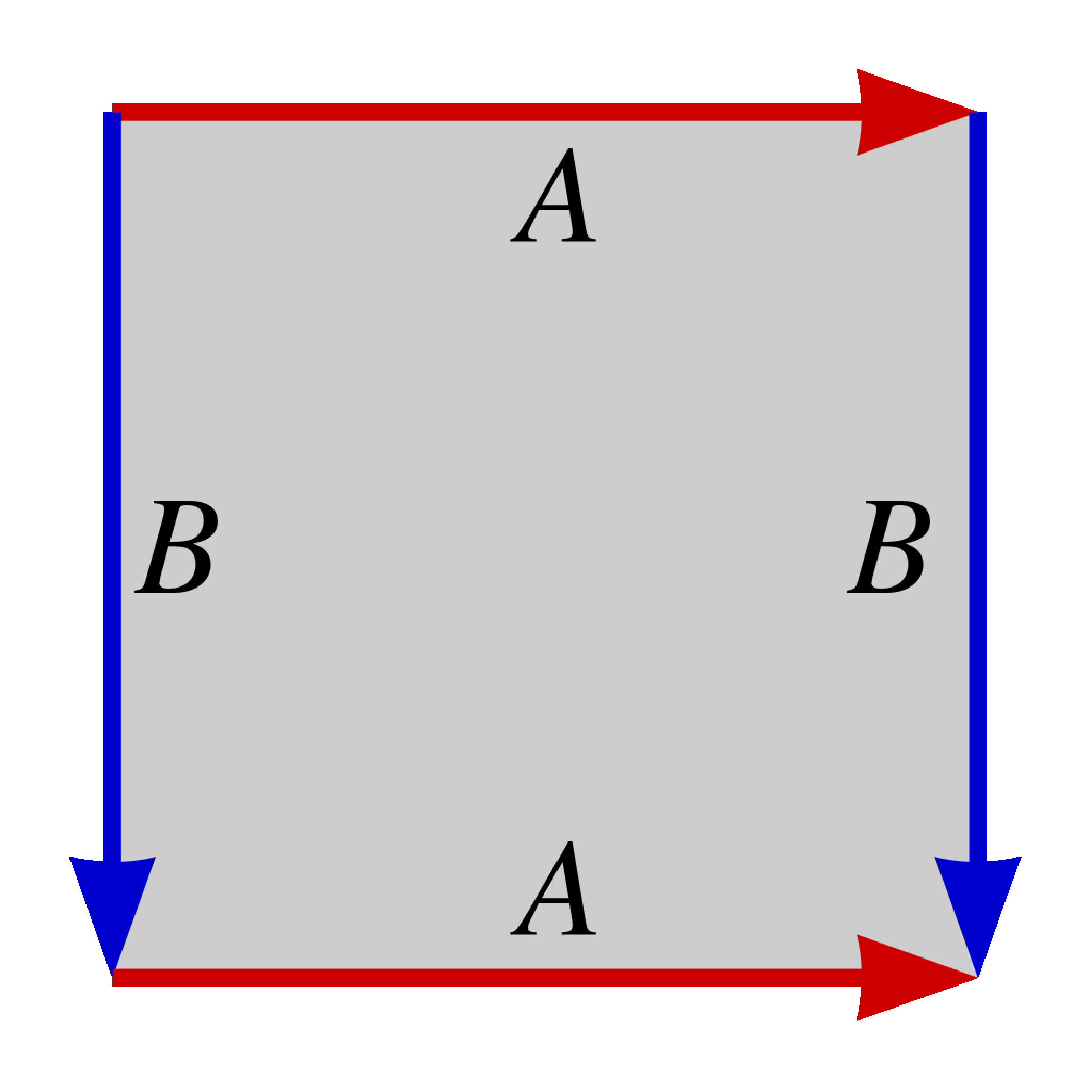}
    \caption{A torus is formed from a square/rectangle by gluing opposite sides, which makes the rectangle a fundamental polygon of the torus. Image Courtesy: Ilmari Karonen, Wikimedia Commons}
\end{figure}

Simply speaking, a Torus in $\mathbb{R}^3$ is just the Cartesian product of two circles, given by $\mathbb{T}^2 = \mathbb{S}^1\times \mathbb{S}^1$ (\cite{hatcher2005algebraic}, pg. 5). In our case, we can form one from a rectangular $m \times n$ board $X$ by \textit{wrapping around} the board along the two dimensions.\footnote{In a physicist's language, this is called following a Periodic Boundary Condition (PBC) in two dimensions.} We must first define an appropriate adjacency matrix $A_{m, n}^\circ$ which endows $X$ with the correct neighborhood structure. Following that, we set $T_{m, n}^\circ = A_{m, n}^\circ - \mathbb{I}_{mn}$ and examine the solution space $\ker(T_{m, n}^\circ)$.

\subsection{Solutions on a Torus}

\begin{prp}\label{prp:kings_graph_torus}
    Define $B_n^\circ \in M_n(\mathbb{Z})$ with $B_{ij}^\circ = 1$ when $i - j \in \{-1, 0, 1\} \pmod{n}$ and $B_{ij}^\circ = 0$ otherwise. Then, \[
        A_{m, n}^\circ = B_m^\circ\otimes B_n^\circ -  \mathbb{I}_{mn}, \qquad
        T_{m, n}^\circ = B_m^\circ\otimes B_n^\circ - 2\mathbb{I}_{mn}.
    \]
\end{prp}
\begin{proof}
    Note that $B_n^\circ$ is the adjacency matrix of the graph $G_n$ from Proposition~\ref{prp:kings_graph} with the extra edge $(n, 1)$. Proceeding in the same manner, we obtain the required adjacency matrix of the toroidal King's Graph $G_{m, n}^\circ$.
\end{proof}

\begin{fct}[Kulkarni, Schmidt and Tsui \cite{kulkarni1999eigenvalues}]
    The eigenvalues of $B_n^\circ$ are $\lambda_k^\circ = 1 + 2\cos(2 k \pi / n)$ for $k = 1, \dots, n$.
\end{fct}
Note that $B_n^\circ$ is a tridiagonal circulant matrix \cite{davis1979circulant} with all non-zero entries being unity, wherefore its eigenvalues are solutions to its associated polynomial
$f(x) = 1+x+x^{n-1}$, which gives the above result.

\begin{prp}\label{prp:eigenvalue_eq_torus}
    The space $\ker(T_{m, n}^\circ)$ is non-trivial if and only if there exist $p, q \in \mathbb{N}$ such that $1 \leq p \leq m$, $1 \leq q \leq n$ and \[
        \left(1 + 2\cos\left(\frac{2p\pi}{m}\right)\right) \left(1 + 2\cos\left(\frac{2q\pi}{n}\right)\right) = 2.
    \]
\end{prp}

\begin{thm}\label{torus}
    Non-trivial $m\times n$ toroidal chessboards satisfying the neighbor-sum property exist if and only if $4 \mid m$ and $6 \mid n$.
\end{thm}
\begin{proof}
    The results follow from Proposition~\ref{prp:eigenvalue_eq_torus} and Theorem~\ref{thm:rational_solutions}. Note that the edge cases where $p = m/2$ or $p = m$ can be eliminated by hand.
\end{proof}

This concludes our discussion on finite generalizations of square chessboards. We showed that for $m\times n$ rectangles, the if and only if condition for the existence of solutions was $2 \mid (m+1)$ and $3\mid (n+1)$. For $m\times n$ torii, solutions exist if and only if $4\mid m$ and $6\mid n$.

There are certain other variations of boards (graphs) where we can consider this problem. This can include infinite boards, higher dimensional analogs or custom connected boards. But before we approach those problems, let us try to switch up the nature of the problem by changing the concept of a neighborhood as defined in \textbf{Section 1}.

\section{The \textit{Neumann Neighborhood}.}
Working with a spatial lattice yields some nice analogs, such as torii. Another way to generate analogs is by redefining the idea of neighbors. We will take a look at the Neumann neighborhood \cite{von1966theory} in this section.
An $m \times n$ chessboard for $m,n \geq 3$, with each square bearing an integer, is said to be a solution to the \emph{Neumann-neighbor-sum (Nns)} problem if each number is the sum of the numbers on the neighboring squares\textemdash where two squares are neighbors if and only if they share a common edge. \textit{Trivial solutions} are defined as before.

Are there non-trivial boards with \textit{Nns} property?

\subsection{On Existence.}
\begin{prp}\label{neumann}
    The set of all $m\times n$ chessboards with the Neumann-neighbor-sum property is precisely $\ker(T_{m, n}^+)$, where \[
        T_{m, n}^+ = B_m\oplus B_n - 3\mathbb{I}_{mn} = B_m\otimes \mathbb{I}_n + \mathbb{I}_m \otimes B_n - 3\mathbb{I}_{mn}.
    \]
\end{prp}
\begin{proof}
    Let $G_{m, n}^+$ be the graph corresponding to an $m\times n$ chessboard with the Neumann neighborhood structure. We claim that its adjacency matrix is \[
        A_{m, n}^+ = B_m\otimes \mathbb{I}_n + \mathbb{I}_m \otimes B_n - 2\mathbb{I}_{mn}.
    \] Indeed, in the graph corresponding to the adjacency matrix $B_m\otimes \mathbb{I}_n$, each square is connected to itself and the squares above and below it. Similarly, in the graph corresponding to $\mathbb{I}_m\otimes B_n$, each square is connected to itself and the squares to its left and right. Adding these two graphs gives $G_{m, n}^+$ except with two self-loops for each square, which we remove by subtracting $2\mathbb{I}_{mn}$ to retrieve the above expression for $A_{m, n}^+$.
\end{proof}

\begin{fct}[Henderson, Pukelsheim and Searle \cite{henderson1983history}]
    The eigenvalues of $A\oplus B$ are $\{\lambda_i + \mu_j\}$, where $\{\lambda_i\}$ are the eigenvalues of $A$, and $\{\mu_j\}$ are the eigenvalues of $B$.
\end{fct}

This gives us an eigenvalue equation unlike those we have seen so far.

\begin{prp}\label{prp:eigenvalue_eq_edge}
    The space $\ker(T_{m, n}^+)$ is non-trivial if and only if there exist $p, q \in \mathbb{N}$ such that \[
        \cos\left(\frac{p\pi}{m+1}\right) + \cos\left(\frac{q\pi}{n+1}\right) = \frac{1}{2}.
    \]
    with $1 \leq p \leq m$ and $1 \leq q \leq n$.
\end{prp}

Now we discuss a particular case $ m=n $. We want solutions of $$ 2\cos\left(\frac{a\pi}{n+1}\right)+2\cos\left(\frac{b\pi}{n+1}\right) = 1$$ with $ a,b\in\{1,2,3,\ldots,n\} $.

\begin{thm}[Conway and Jones \cite{conway1976trigonometric}]\label{coscomb}
    Suppose we have at most four distinct rational multiples of $\pi$ lying strictly between 0 and $\pi / 2$ for which some rational linear combination of their cosines is rational but no proper subset has this property. Then the appropriate linear combination is proportional to one from the following list:
 {\setlength{\columnsep}{-1.1cm}
 \begin{multicols}{2}
    \begin{small}

    \begin{itemize}
    \renewcommand\labelitemi{} 
	\setlength{\leftmargin}{0cm}
        \item  \hspace*{-1em} $\cos (\frac{\pi}{3})=\frac{1}{2}$
		\item \hspace*{-1em} $\cos (\varphi)-\cos (\frac{\pi}{3} -\varphi)-\cos (\frac{\pi}{3} +\varphi)=0$,\\ where $0<\varphi<\frac{\pi}{6}$
		\item \hspace*{-1em} $\cos (\frac{\pi}{5})-\cos (\frac{2\pi}{5})=\frac{1}{2}$
		\item \hspace*{-1em} $\cos (\frac{\pi}{7})-\cos (\frac{2\pi}{7})+\cos (\frac{3\pi}{7})=\frac{1}{2}$
		\item \hspace*{-1em} $\cos (\frac{\pi}{5})-\cos (\frac{\pi}{15})+\cos (\frac{4\pi}{15})=\frac{1}{2}$
		\item \hspace*{-1em} $\cos (\frac{2\pi}{5})-\cos (\frac{2\pi}{15})+\cos (\frac{7\pi}{15})=-\frac{1}{5}$ 
		\item \hspace*{-1em} $\cos (\frac{\pi}{7})+\cos (\frac{3\pi}{7})-\cos (\frac{\pi}{21})+\cos (\frac{8 \pi}{21})=\frac{1}{2}$
		\item \hspace*{-1em} $\cos( \frac{\pi}{7})-\cos (\frac{2\pi}{7})+\cos (\frac{2\pi}{21})-\cos( \frac{5\pi}{21})=\frac{1}{2}$
		\item \hspace*{-1em} $\cos (\frac{2\pi}{7})-\cos (\frac{3\pi}{7})-\cos (\frac{4\pi}{21})-\cos (\frac{10\pi}{21})=-\frac{1}{2}$
		\item \hspace*{-1em} $\cos (\frac{\pi}{15})-\cos( \frac{2\pi}{15})-\cos (\frac{4\pi}{15})+\cos (\frac{7\pi}{15})=-\frac{1}{2}$ 
	\end{itemize}
        
    \end{small}
	\end{multicols}}
\end{thm}

	\textbf{Note:} If $ a\ge \frac{n+1}{2} $, $ a':=n-a\le\frac{n+1}{2} \implies \cos\left(\dfrac{a\pi}{n+1}\right) = -\cos\left(\dfrac{a'\pi}{n+1}\right).$

\begin{thm}
	The equation has a solution if and only if $ 5\mid n+1 $ or $ 6 \mid n+1 $.
\end{thm}

\begin{proof}
	By the Note, the problem reduces to finding solutions of $$ \pm 2\cos\left(\frac{a'\pi}{n+1}\right) \pm 2\cos\left(\frac{b'\pi}{n+1}\right) = 1$$ with $ 0<a',b'\le \frac{n+1}{2} $.
	
	Also, if $ a'=\frac{n+1}{2}$, $ 2\mid n+1 $ and this implies $$ \cos\left(\frac{b\pi}{n+1}\right) = \frac{1}{2} \implies \frac{b'}{n+1}= \frac{1}{3} \implies 2,3\mid n+1$$ 
    implying $6\mid n+1$.
    
    The same analysis holds if one of the term is zero i.e,, $$ 2\cos\left(\frac{a'\pi}{n+1}\right) = 0  \implies a' = \frac{n+1}{2}$$
    and hence $6\mid n+1$.
	
	So, we can assume that none of the terms is zero and $ 0<\frac{a'}{n+1},\frac{b'}{n+1}<\frac{1}{2} $. By applying Theorem~\ref{coscomb} the only 2 term relation is $$ \cos\left(\frac{\pi}{5}\right) - \cos\left(\frac{2\pi}{5}\right) = \frac{1}{2} \implies 2\cos\left(\frac{\pi}{5}\right) + 2\cos\left(\frac{3\pi}{5}\right) = 1 $$
	which implies $ 5\mid n+1 $. So, if the solution of the equation exists, then either $ 6 \mid n+1 $ or $ 5\mid n+1 $.
	
	For the converse,
	if $ n+1=6k $, we have $ 2\cos\left(\frac{3k\pi}{n+1}\right) + 2\cos\left(\frac{2k\pi}{n+1}\right) = 1 $.\\
 If $ n+1=5k $, we have $ 2\cos\left(\frac{k\pi}{n+1}\right) + 2\cos\left(\frac{3k\pi}{n+1}\right) = 1 $
\end{proof}

We can do the same analysis for the $m \ne n$ case and get that the equation has a solution if and only if $5$ divides both $m+1$ and $n+1$ or $3\mid (m+1)$ and $2\mid (n+1)$ or vice versa.

Further, due to the transformation $T_{n,n}^{+}$ being symmetric and having zero eigenvalue with multiplicity $2$ (whenever solutions exist), the kernel is two-dimensional whenever it is non-trivial.

\subsection{Solutions of Nns.}

We will take a look at solutions for $n=4$, which is the smallest value of $n\geq 3$ such that $5$ divides $n+1$.

\[\def\arraystretch{1.7}\arraycolsep=1pt
        \symking^{(1)+}_{4\times 4} =
        \begin{array}{|c|c|c|c|}\hline
            \phantom{-}0\phantom{-} & 1 & 1 & \phantom{-}0\phantom{-} \\\hline
            -1 & 0 & 0 & -1\\\hline
            -1 & \phantom{-}0\phantom{-} & \phantom{-}0\phantom{-} & -1\\\hline
             0 & 1 & 1 & 0 \\\hline
        \end{array}
        \qquad \text{and}\qquad
        \def\arraystretch{1.7}\arraycolsep=1pt
        \symking^{(2)+}_{4\times 4} =
        \begin{array}{|c|c|c|c|c|}\hline
            \phantom{-}1\phantom{-} & \phantom{-}0\phantom{-} & \phantom{-}0\phantom{-} & \phantom{-}1\phantom{-} \\\hline
            1 & -1 & -1 & 1 \\\hline
            1 & -1 & -1 & 1 \\\hline
            1 & 0 & 0 & 1 \\\hline
        \end{array}
        \; .
    \]

Note that unlike the original case (i.e,, the Moore neighborhood), here the elements of the basis are not transposes of each other. The $n=4$ solutions can be mirrored to get extended solutions for larger boards.

\begin{figure}[htp!]
    \centering
    \includegraphics[width=1.0\textwidth]{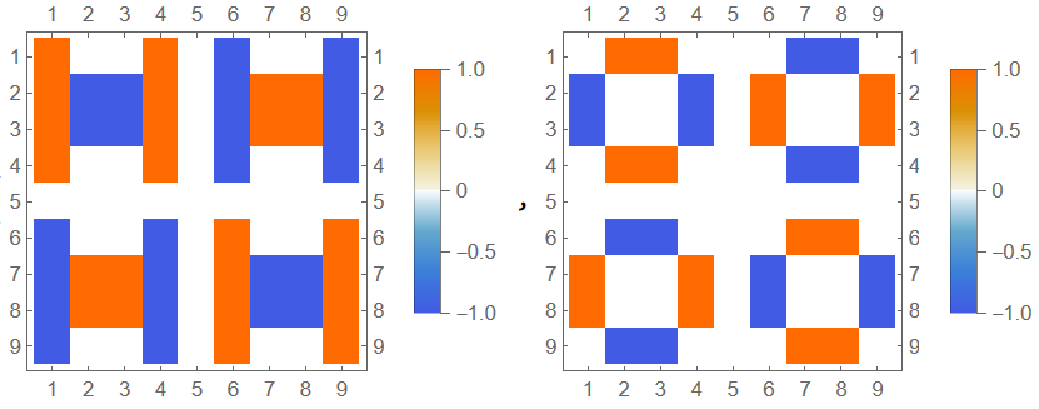}
    \caption{Solutions for $n=9$. The color bar shows that the blue regions are $-1$, the orange ones are $1$, and the rest are zeroes. Notice how the $n=4$ solutions are neatly extended to form this basis.}
\end{figure}

For the case when $6$ divides $n+1$, check that the solutions for the Moore neighborhood also work for the Neumann neighborhood and can be similarly extended from the $n=5$ case.

\subsection{An interesting problem on harmonic functions.}
A \textit{discrete harmonic function} \cite{heilbronn1949discrete} on a graph $G = (\mathcal{V}, \mathcal{E})$ is defined as follows:

\begin{defi}\label{harmonic}
    A function $f:\mathcal{V}\to \mathbb{R}$ is harmonic at a node/vertex $x\in\mathcal{V}$ if it satisfies the following relation
    \[f(x) = \frac{\sum_{\{x,y\}\in\mathcal{E}}f(y)}{deg(x)}\]
    where $deg(x)$ is the degree of the vertex $x$.
\end{defi}

Consider a toroidal chessboard with the Neumann neighborhood condition i.e,, two squares are neighbors if and only if they share a common edge. To define the associated graph $G$, we identify the squares with vertices and draw an edge between every neighbor of the torus. This creates a graph where every vertex has degree $4$.

Since every vertex is a representation of a cell on the toroidal board, we can denote a pair of coordinates $(p,q)$ to represent it, where $1\leq p\leq m,\;1\leq q\leq n$. It is easy to then find an appropriate transformation $T^{\circ+}_{m,n}$ whose kernel contains solutions of the modified neighbor sum equation
\[f(p,q+1)+f(p,q-1)+f(p+1,q)+f(p-1,q) = 4f(p,q)\]
where if $(p,q)$ represents vertex $x$, then the tuples $(p\pm 1, q)$, $(p,q\pm 1)$ represents its Neumann neighbors. Note that this is not a neighbor sum problem, rather the mean tuple is one-fourth of the neighbor sum. Call this a \textit{Neumann-neighbor-average} property of the vertex $x = (p,q)$.

\begin{prp}\label{prp:harmonic-transform}
    The set of all $m\times n$ toroidal chessboards with the Neumann-neighbor-average property is precisely $\ker(T^{\circ+}_{m,n})$, where
    \[T^{\circ+}_{m,n} = B^{\circ}_m\oplus B^{\circ}_n-6\mathbb{I}_{m,n}.\]
\end{prp}
\begin{proof}
    The proof is along the same lines as that of Proposition~\ref{neumann}. Here the matrices $B^{\circ}_m$, $B^{\circ}_n$ correspond to adjacency matrices of circular graphs with self-loops. Further, in the final transformation $T^{\circ+}_{m,n}$, we require the diagonal elements to be $(-4)$ keeping all other elements constant, so we subtract $6\mathbb{I}_{m,n}$.
\end{proof}

\begin{prp}
    The space $\ker(T_{m, n}^{\circ+})$ is non-trivial if and only if there exist $p, q \in \mathbb{N}$ such that $1 \leq p \leq m$, $1 \leq q \leq n$ and \[
        \cos\left(\frac{2p\pi}{m}\right) + \cos\left(\frac{2q\pi}{n}\right) = 2.
        \]
\end{prp}

This eigenvalue equation has the trivial solution $\cos\left(\frac{2p\pi}{m}\right) = \cos\left(\frac{2q\pi}{n}\right) = 1$, which gives $p=m,\; q=n$. So, the kernel is one-dimensional with the corresponding eigenvalues of $B^{\circ}_m$ and $B^{\circ}_n$ are both $3$.

\begin{fct}
    The eigenvectors of $A\oplus B$ are $\{x_i\otimes y_j\}$ where $\{x_i\}$ are the eigenvectors of $A$ and $\{y_j\}$ are the eigenvectors of $B$.
\end{fct}

Note that eigenvectors of $B_m^{\circ}$ over $\mathbb{Z}^m$ for the eigenvalue $3$ are of the form $c\mathbbm{1}_m$, where $c\in\mathbb{Z}\backslash\{0\}$ is a constant. Then the corresponding unique solution (in vectorized format) is $c'\mathbbm{1}_m\otimes \mathbbm{1}_n = c'\mathbbm{1}_{mn},\;\; c'\in\mathbb{Z}\backslash\{0\}$, which is just a constant.

This gives a neat result as mentioned below.
\begin{thm}
    Discrete harmonic functions on a toroidal graph are constant functions.
\end{thm}

\textit{Remark.} If the same problem is specified on a finite square lattice, we could imagine the presence of a phantom boundary such that the mean-value property as defined is consistent. As the harmonic function takes the value zero on the boundary, by the maximum modulus principle for discrete harmonic functions \cite{heilbronn1949discrete}, the only value it can take in the interior is identically zero. This can be shown with our linear algebra machinery and is a nice connection between the two ideas.

\section{Board extends to Infinity.}
At this juncture, it seems very natural to ask the same question for infinite chessboards. However, since an infinite board has fewer restrictions, it is much easier to answer this question.

A semi-infinite chessboard is one where the board is infinite in only one direction along the $x$ and $y$ axes. In that case, every square on the board can be expressed as a tuple $(i,j), \; i,j\in\mathbb{N}$. We do this numbering akin to that of rows and columns from the cell $(1,1)$ with only three neighbors sharing a common vertex or edge.

An infinite chessboard, which is infinite in both directions along the x and y axes is numbered as follows. Take a row and column and call them the 0-th row and column. Columns left of the 0-th column (respectively rows above the 0-th row) will be numbered by the negative integers while those right of (respectively below) would be numbered by the positive integers. The element at the intersection of the $i$-th row and $j$-th column is represented by the tuple $(i,j)$.

 \subsection{Semi-infinite chessboard}
 For a semi-infinite board, denote the value at any cell enumerated by $(i,j)$ as $x_{ij}$.
 
 We can take any two sequences $\{a_i\}_{i=1}^\infty$ and $\{b_i\}_{i=1}^\infty$ (with $a_1=b_1$) to fill up the positions $\{(1,i)\}_{i=1}^\infty$ and $\{(i,1)\}_{i=1}^\infty$ respectively (enumeration begins from top left). Notice that the values at $(1,1), (1,2)$, and $(2,1)$ fix $(2,2)$, call this value $x_{22}$.

 $x_{22}$ along with the other given values fixes $x_{23}$ and $x_{32}$, and it is easy to check that done recursively, this fixes row 2 and column 2.

 Note that in the first scenario, we could have used a phantom boundary beyond the first row and column to show that the existence of two adjacent filled rows and columns fixes the third row and column for semi-infinitely long axes. By recursively fixing elements, we can generate the solution. We are using a matrix notation just for clarity (this is not a matrix!).

 \[\symking_{\mathbb{N}} =
        \begin{pmatrix}
            2 & 3 & 5 & 7 & 11 & 13 &.&.&.\\
            3 & -4 & -3 & 2 & -8 & -3 &.&.&.\\
            5 & -3 & -16 & 3 & 3 & -42 &.&.&.\\
            8 & 1 & 3 & -15 & 37 & 29 &.&.&.\\
            13 & -8 & -1 & 42 & -86 & 99 &.&.&.\\
            21 & -7 & -46 & 29 & 121 & -428 &.&.&.\\
            .&.&.&.&.&.&.&.\\
            .&.&.&.&.&.&.&.\\
            .&.&.&.&.&.&.&.\\
        \end{pmatrix}
        \]

 \subsection{Infinite chessboard}

For an infinite board, it is again a matter of choosing the sequences $\{a_i\}_{i\in \mathbb Z\setminus \{0\}}$, $\{b_i\}_{i\in \mathbb Z\setminus \{0\}}$, $\{c_i\}_{i\in \mathbb Z\setminus \{0\}}$, and $\{d_i\}_{i\in \mathbb Z\setminus \{0\}}$ satisfying
\begin{align*}
    a_1=b_1,\qquad a_{-1}=d_1,\qquad c_1=b_{-1},\qquad c_{-1}=d_{-1}
\end{align*}
and putting them along any two rows and columns as shown in the figure below.

\[
\symking_{\mathbb Z}=
    \begin{pmatrix}
        .&.&.&.&\astrosun &\astrosun &.&.&.&.\\
        .&.&.&.&\astrosun &\astrosun &.&.&.&.\\
        .&.&.&.&\astrosun &\astrosun &.&.&.&.\\
        .&.&.&.&\astrosun &\astrosun &.&.&.&.\\
        \astrosun &\astrosun &\astrosun &\astrosun &\astrosun &\astrosun &\astrosun &\astrosun &\astrosun &\astrosun\\
        \astrosun &\astrosun &\astrosun &\astrosun &\astrosun &\astrosun &\astrosun &\astrosun &\astrosun &\astrosun\\
        .&.&.&.&\astrosun &\astrosun &.&.&.&.\\
        .&.&.&.&\astrosun &\astrosun &.&.&.&.\\
        .&.&.&.&\astrosun &\astrosun &.&.&.&.\\
        .&.&.&.&\astrosun &\astrosun &.&.&.&.\\
    \end{pmatrix}
\]

If the entries marked with $\astrosun$ are filled up with any sequences, they would divide the infinite board into four disjoint semi-infinite boards (disjoint since any entry in each of the four quadrants would only depend upon the already filled numbers and not entries in any other quadrant). Further, each semi-infinite board contains a unique solution if two adjacent rows and columns are specified, which is the case here. Notice that the entries marked with $\astrosun$ can be first filled up by any numbers without already violating the neighbor-sum property.

\section{In Conclusion.}

Starting from a Math Olympiad problem, we generated a variant termed the neighbor sum problem, where we look for a non-constant (\textit{a.k.a} non-trivial) vector of integers represented on a chessboard, such that the value in each square is equal to the sum of its neighbors.

Using some linear algebra machinery we arrived at a trigonometric diophantine equation, solving for which 2D finite boards required some p-adic analysis. This provided necessary and sufficient conditions for non-trivial solutions to exist on square, rectangular or toroidal boards.

Eventually, we checked out the case of infinite boards, where generating solutions is quite easy, as well as finite boards under the Neumann neighborhood\textemdash which required a result by Conway and Jones \cite{conway1976trigonometric}.

Nevertheless, it appears imperative to us to set a grand generalization of this problem, whose specific cases would be the ones we have treated in this paper. In the following two sections, we will look at some rigorous mathematical insight into solving the relevant eigenvalue equations on higher dimensional grids like hypercubes. At the end, we pose \textbf{Question 31}, which provides an avenue for further research.

\section{A brief Digression\textemdash some more Algebra.}
    We wanted to get the integer $2$ as products of algebraic integers of the form
    $$1 + 2\cos\left(\frac{2\pi p}{2(n+1)}\right) = 1 + \zeta_m^a + \zeta_m^{-a}=:\lambda_{a,m}$$
    where $\frac{a}{m}$ is the reduced form of $\frac{p}{2(n+1)}$ . Naturally, we hoped that if one of these algebraic integers involving a primitive $m$-th root of unity appeared in the product, all of its conjugates would too. This is true in two dimensions, but a counter-example in three dimensions is as follows:
    \[\left(1 + 2\cos\left(\frac{2\pi}{24}\right)\right)\left(1 + 2\cos\left(\frac{22\pi}{24}\right)\right)\left(1 + 2\cos\left(\frac{20\pi}{24}\right)\right) = 2.\]
    So, to get a sufficient condition for the existence of a solution we calculated the product of the conjugates i.e, the usual field norm \cite{rotman2010advanced} of $\lambda_{a,m}$ say $g(m)$, which would be an integer. If the product of some of these norms equals 2 and the total number of terms counting conjugates equals $d$, we get a solution.
    
    Motivated by this, let us formally define
    $$ g(m):= \displaystyle\prod_{\begin{array}{c} 1 \leq a < \frac{m}{2}\\ \gcd(a, m) = 1 \end{array}} \left(1 + 2\cos\left(\frac{2\pi a}{m}\right)\right).$$
	
    Observe that if we can write $ 2 $ as the product of $ g(m_i) $'s where each $ m_i \mid 2(n+1)$ and the ``total length" of the product is $d$, we will have a solution in $d$ dimensions,
    $$\displaystyle\prod_{i=1}^d \left(1+2\cos \left(\frac{p_i \pi}{n+1}\right)\right)=2.$$
    Here, by the length of $ g(m_i) $, we mean the number of terms appearing in the product, i.e, $ \dfrac{\phi(m_i)}{2} $, and by total length we mean $ \displaystyle\sum_{m_i} \frac{\phi(m_i)}{2} $. As we'll see in Theorem \ref{necessary}, $ 3\mid (n+1) $ is a necessary condition, if $2\mid (n+1)\implies 4\mid 2(n+1)$, we can choose $ m_1=6 $ and $ m_2=\cdots=m_d=4 $. This gives us a solution. To generalize the idea and get a better sufficient condition, we calculate $ g(m) $ for different primes.

    The following theorem provides an expression of $g(m)$ in terms of the minimal polynomial \cite{lehmer1933note,watkins1993minimal} of $2\cos\left(\frac{2\pi}{m}\right)$ and the $ m $-th cyclotomic polynomial (\cite{isaacs2009algebra} pg. 308).
	
	\begin{thm}\label{g}
        The function $g$ satisfies
		$$ g(m) = \Psi_m(-1)(-1)^{\frac{\phi(m)}{2}}=\Phi_m(\zeta_3)(\zeta_3)^{\frac{-\phi(m)}{2}}(-1)^{\frac{\phi(m)}{2}}  $$ 
		where $ \Psi_m $ is the minimal polynomial of $ 2\cos\left(\frac{2\pi}{m}\right)$, $\Phi_m $ is the $ m $-th cyclotomic polynomial, and $ \zeta_3 $ is a primitive 3rd root of unity.
	\end{thm}  

	\begin{proof}
        Taking
		$$ \Psi_m(x) = \displaystyle\prod_{\begin{array}{c} 1 \leq a < \frac{m}{2}\\ \gcd(a, m) = 1 \end{array} } \left(x - 2\cos\left(\frac{2\pi a}{m}\right)\right).$$
		and putting $ x=-1 $, we get the first equality.
        
        It is easy to check that
        $$ \Psi_n\left(z + z^{-1}\right) = z^{-\frac{\phi(n)}{2}}\Phi_n(z)$$ 
        and hence, putting $ z=\zeta_3 $ we get the second equality.
	\end{proof}
	
	Now using this we can explicitly calculate $ g(m) $ for any $ m\in\mathbb{N}, m>3 $. We'll show some of the relevant calculations here.
	\begin{thm}\label{gp}
		Let $ p>3 $ be a prime. Then,
            $$g(p)=\left(\frac 3p\right)=\begin{cases}
                1 & \text { if } p\equiv     \pm 1 \pmod{12} \\
				-1 & \text { if } p\equiv \pm 5 \pmod{12}
            \end{cases}$$
		where $\left(\frac{\textcolor{white}{\cdot}}{\textcolor{white}{\cdot}}\right)$ is the Legendre symbol. Moreover, $g(2p)=1$ if $p\equiv 5 \pmod{12}$.
	\end{thm}
	\begin{proof}
		$ \Phi_p(x) = 1 + x + \cdots + x^{p-1} $ and $ \Phi_{2p}(x) = 1 - x + \cdots + x^{p-1} $.\\
		$$ (-1)^{\frac{p-1}{2}}=\left\{\begin{array}{cc}
			1 & \text { if } p\equiv  1 \pmod{4} \\
			-1 & \text { if } p\equiv -1 \pmod{4}
		\end{array}\right.$$. $$ {\zeta_3}^{-\left(\frac{p-1}{2}\right)}\; \Phi_p(\zeta_3)=\left\{\begin{array}{cc}
		1 & \text { if } p\equiv  1 \pmod{3} \\
		-1 & \text { if } p\equiv -1 \pmod{3}
	\end{array}\right. $$.\\
	Also, $ \phi(p)=\phi(2p) $ for any odd prime and $ {\zeta_3}^{-\left(\frac{p-1}{2}\right)}\Phi_{2p}(\zeta_3)= 1 $ if $ p\equiv -1 \pmod{3} $. This along with Theorem~\ref{g} gives us the proof.
	\end{proof}

\section{Higher Dimensions and other generalizations.}
Using the tools developed in the previous sections, we will look at some nice results for higher dimensional analogs of the neighbor sum problem, such as on hypercubes. This is a natural extension of all our work until now, yet the conditions for existence or explicit solutions might be completely different from what we have observed. The problem can also be extended to arbitrary graphs or lattice-subsets in $\mathbb Z^{d}$, but the problem becomes too complicated to even get a flavor of the kinds of solutions expected. Different tools might be needed to tackle them.

For example, start by considering the same problem on a $d$-dimensional hypercube. The analog of the equation in Proposition~\ref{prp:eigenvalue_eq} for this case is
$$\prod_{i=1}^d \left(1+2\cos \left(\frac{p_i \pi}{n+1}\right)\right)=2$$
which is equivalent to trying to find the number of solutions of the equation
$$\sum_{i=1}^d v \left(\alpha_i+1+\alpha_i^{-1}\right)=1$$
which now has a lot more combinations than was possible in the $d=2$ case. A criterion for existence of solutions is no longer as straightforward as on rectangles or squares. However, a necessary condition for existence is easy to figure out.

\begin{thm}\label{necessary}
    If an $n^d$ board satisfies the neighbor sum problem, then $3|(n+1)$.
\end{thm}
\begin{proof}
    As already established, if a $n^d$ board satisfies the problem, then the equation
    $$\sum_{i=1}^d v \left(\alpha_i+1+\alpha_i^{-1}\right)=1$$
    has a solution.

    But, this means that there is a $k$ such that $v(\alpha_i+1+\alpha_i^{-1})>0$, hence completing the proof.
\end{proof}

\begin{thm}\label{3cube}
    There are solutions of
    $$\prod_{i=1}^{d=3} \left(1+2\cos \left(\frac{p_i \pi}{n+1}\right)\right)=2$$
    if and only if $6 \mid (n+1)$ or $ 15\mid (n+1)$. So, an $n\times n\times n$ board has a solution if and only if $6 \mid (n+1)$ or $ 15\mid (n+1)$.
\end{thm}
\begin{proof}
Let $n+1= 6k$ for some $k\in\mathbb{N}$. $p_1 = 2k$, $p_2=p_3=3k$ solves the eigenvalue equation.\\
Let $ n+1 =15k $. Then $p_1 = 5k$, $p_2=3k$, $p_3=9k$ is a solution.

Now to prove the converse, we already had $3\mid (n+1) $ by the previous theorem. Assume $2\nmid (n+1)$. Let the solution of the equation be $$ \left(1+2\cos \left(\frac{p_0 \pi}{n+1}\right)\right)\left(1+2\cos \left(\frac{p_1 \pi}{n+1}\right)\right)\left(1+2\cos \left(\frac{p_2 \pi}{n+1}\right)\right)=2 $$
WLOG, we obtain two possibilities by the 2-adic valuation: $1+0+0$ and $\frac{1}{2}+ \frac{1}{2} +0 $. Since, $v(\alpha_i+1+\alpha_i^{-1})= \frac{1}{2} $ gives us the 12-th roots of unity, we get a contradiction that $2\nmid (n+1)$. So, the only possible case is $1+0+0$. So, assume that $p_0= \frac{n+1}{3}$, i.e, the first term is $2$. Thus, the problem reduces to
\begin{align*}
    &\left(1+2\cos \left(\frac{p_1 \pi}{n+1}\right)\right)\left(1+2\cos \left(\frac{p_2 \pi}{n+1}\right)\right) = 1\\
    \implies & 2\cos\left(\frac{p_1 \pi}{n+1}\right) + 2\cos\left(\frac{p_2 \pi}{n+1}\right) + 4 \cos\left(\frac{p_1 \pi}{n+1}\right)\cos\left(\frac{p_2 \pi}{n+1}\right) = 0\\
    \implies &\cos\left(\frac{p_1 \pi}{n+1}\right) + \cos\left(\frac{p_2 \pi}{n+1}\right) + \cos\left(\frac{p_1+p_2}{n+1}\;\pi\right) + \cos\left(\frac{p_1-p_2}{n+1}\;\pi\right) = 0.
\end{align*}
Note that $2\nmid(n+1)$ implies that no single term in the equation is $0$ and $p_1\ne p_2$. $p_1,p_2\ne 0 $ implies that all four angles in the equation are distinct. By Theorem \ref{coscomb}, we have a description of the possible angles whose cosine sum is 0.

We can reduce the angles to be less than $\pi/2$ by the following:

If $\dfrac{\pi}{2} < \theta < \pi $, $\cos\left(\theta\right) = - \cos\left(\theta'\right)$ where $\theta'= \pi - \theta < \dfrac{\pi}{2}$.

If $\theta> \pi$, $\cos\left(\theta\right)= \cos\left(\pi + \theta - \pi\right) = -\cos\left(\theta-\pi\right)$.

Note that by this we can make the first, second and fourth angle to lie strictly between $0$ and $\pi/2$. For the third angle only possibility where this cannot be done remains:\\
If $p_1 + p_2 = n+1$, the equation gives us $p_1 = p_2$, contradiction.\\

Then we will have the following cases:\\
Case 1: The angles are not distinct $\implies$ they should cancel each other in two-two pair (we can eliminate other possibilities by observing Theorem \ref{coscomb})\\
Case 2: The angles are distinct $\implies$ some subsequence will have rational sum $\implies$ the LHS is difference of some subsequence from the following list:

  \begin{itemize}
	 	\item[a.] $\cos \left(\frac{\pi}{3}\right)=\frac{1}{2}$
		\item[b.] $\cos \left(\frac{\pi}{5}\right)-\cos \left(\frac{2\pi}{5}\right)=\frac{1}{2}$
		\item[c.] $\cos \left(\frac{\pi}{7}\right)-\cos \left(\frac{2\pi}{7}\right)+\cos \left(\frac{3\pi}{7}\right)=\frac{1}{2}$
		\item[d.] $\cos \left(\frac{\pi}{5}\right)-\cos \left(\frac{\pi}{15}\right)+\cos \left(\frac{4\pi}{15}\right)=\frac{1}{2}$.
  \end{itemize}

For the second case, we see that the only possible case where the terms don't cancel each other is $d-a$ which implies $5\mid(n+1)$. Note that the $c-a$ case is not possible because $\frac{\pi}{3}$ cannot be written as sum of other three angles.

For the case that the pair of two terms cancel each other, we have several subcases:

\textbf{Case 1: (First-Second, Third-Fourth Pair). }\\
$$\cos\left(\dfrac{p_1 \pi}{n+1}\right) =- \cos\left(\dfrac{p_2 \pi}{n+1}\right)$$ and $$\cos\left(\dfrac{p_1+p_2}{n+1}\;\pi\right) =-\cos\left(\dfrac{p_1-p_2}{n+1}\;\pi\right).$$\\
The first equation implies $$p_1 + p_2 = n+1 \implies \cos\left( \dfrac{p_1 - p_2}{n+1}\;\pi\right) = 1$$ which implies $p_1=p_2$ contradiction.

\textbf{Case 2: (First-Third, Second-Fourth Pair).}\\
$$\cos\left(\dfrac{p_1 \pi}{n+1}\right) =-\cos\left(\dfrac{p_1+p_2}{n+1}\;\pi\right)\implies n+1 = p_1+p_2 \pm p_1$$
$p_2\ne n+1 \implies 2p_1+p_2 = n+1\implies 4p_1 + 2p_2 = 2(n+1)$ . 
$$\cos\left(\dfrac{p_2 \pi}{n+1}\right) =-\cos\left(\dfrac{p_1-p_2}{n+1}\;\pi\right)\implies n+1 = \pm(p_1 - p_2) + p_2$$
$p_1\ne n+1 \implies 2p_2-p_1 = n+1\implies 5 p_1 = n+1 \implies 15\mid (n+1)$.

\textbf{Case 3: (First-Fourth, Second-Third Pair).}\\
This case is similar to the previous case by changing the roles of $p_1 $ and $p_2$.
\end{proof}

It is also not very difficult to arrive at a sufficient condition for the problem for $n^d$ board. We have discussed the required language in Section 8.

By using Theorem \ref{gp}, we give the sufficient condition for the existence of solutions.

\begin{thm}
   Let $ n+1 = 3^{a_0} \; p_1^{a_1}\cdots p_r^{a_r} \; q_1^{b_1}\cdots q_s^{b_s}$ where $ q_i\equiv 7 \pmod{12} $ and $ p_i \not\equiv 7 \pmod{12}$. Then
   \begin{itemize}
       \item[1.] If any $p_i =2$, a solution of the form $g(6)g(4)\cdots g(4)=2$ exists.
       \item[2.] If all $p_i$'s are odd primes, then if there are integers $x_1,\cdots,x_{r},y_1,\cdots,y_{s}\ge 0$ such that 
$$ x_1 \left (\frac{p_{1}-1}{2}\right ) + \cdots + x_{r}\left (\frac{p_r-1}{2}\right ) + 2 y_1 \left (\frac{q_1-1}{2}\right ) + \cdots + 2 y_{s}\left (\frac{q_{s}-1}{2}\right ) = d-1, $$
then a solution exists and is given by
$$g(6)\; g(p'_1)^{x_1}\cdots g(p'_r)^{x_r} \; g(q_1)^{2y_1}\cdots g(q_s)^{2y_s}=2$$
where $ p'_i =2p_i$ if $ p_i \equiv 5 \pmod{12} $ and $ p'_i = p_i$ otherwise.
   \end{itemize}
\end{thm}

\begin{proof}
    It is enough to explicitly compute the form of the solutions.
\end{proof}

We have also previously established that for $d=2$, whenever we have a solution, the solution space is two dimensional. This motivates us to ask the question for higher dimensional analogs of the problem. For this question, we do not have any useful results to present. We give the sequence $\{a_n^d\}_{n\ge 2}$ of numbers for solutions of $n^d$ boards obtained numerically.

$d=3$ :
\\0, 0, 0, 3, 0, 0, 0, 0, 0, 15, 0, 0, 6, 0, 0, 3, 0, 0, 0, 0, 0, 15, 0, 0, 0, 0, 0, 9, 0, 0, 0, 0, 0, 15, 0, 0, 0, 0, 0, 3, 0, 0, 6, 0, 0, 15, $\ldots$

$d=4$ :
\\0, 0, 0, 4, 0, 0, 0, 0, 0, 88, 0, 0, 24, 0, 0, 4, 0, 0, 0, 0, 0, 136, 0, 0, 0, 0, 0, 220, 0, 0, 0, 0, 0, 88, 0, 0, 48, 0, 0, 52, 0, 0, 24, 0, 0, 136, $\ldots$

$d=5$ :
\\0, 0, 0, 5, 0, 0, 0, 0, 0, 335, 0, 0, 480, 0, 0, 485, 0, 0, 540, 0, 0, 1295, 0, 0, 0, 0, 0, 1865, 0, 0, 0, 0, 0, 815, 0, 0, 0, 0, 0, 1385, 0, 0, 480, 0, 0, 2255, $\ldots$


While these sequences show some patterns, it is hard to find one for general $d$. The articles \cite{roots, rational} provide some techniques and upper bounds of the order of the roots of unity appearing in the eigenvalue equation. Thus, one can try to use the bounds to calculate the possible solution in finite cases with the help of a computer and find generalized results.

In our paper, we defined a transformation whose kernel contains the relevant solutions. The key step was to decompose the graph (subset of $\mathbb{Z}^2$) into Cartesian products of line-graphs (or cycle graphs), which translates neatly over to the adjacency matrix representations through Kronecker products and sums. Similar decompositions for connected graphs can help in reducing the difficulty of the problem and introduce several interesting ideas. One might need a thorough understanding of spectral graph theory and linear algebra to tackle such general problems.

For the sake of completion, we conclude this section by posing a generalized neighbor sum problem on connected graphs.\footnote{It is enough to study it for connected graphs since a disconnected graph satisfies this property if and only if all its connected components do.}
\begin{q}
    For a connected graph $\mathcal G=(\mathcal V,\mathcal E)$ and for some abelian group $G$, call a function $f : \mathcal V \to G$ to satisfy the neighbor sum property if
    \[f(x) = \sum_{\{x,y\}\in\mathcal{E}}f(y)\]
    for all $x\in \mathcal V$. Classify all $\left (\mathcal G, G, f\not \equiv 0\right)$ satisfying this property.
\end{q}

\section{Acknowledgments.}
We would like to acknowledge the anonymous reviewers for providing some useful references and in particular pointing out that the technique of treating similar problems on graphs by studying the eigenvalues of adjacency matrices appears in the classic \textit{Proofs from THE BOOK} \cite{book_proof} in the chapter \textit{Communicating without Errors}-in fact it also appears in other texts like \cite{algcomb} and \cite{titu}. We would also like to thank Dr. David E. Speyer for several meaningful discussions. We would also like to thank our batch mates in CS for helping out with the generation of numerical solutions and relevant sequences. We are indebted to their help in successfully completing this paper early.

\bibliographystyle{plain}
\bibliography{Chotushkone}
\end{document}